\documentclass[12pt,a4paper]{article}
\usepackage[T1]{fontenc}
\usepackage{amsmath}
\usepackage{amssymb}
\usepackage{textcomp}
\usepackage[english,french]{babel}
\usepackage[all,dvips]{xy}
\usepackage[applemac]{inputenc}

\newtheorem{thm}{Theorem}[section]
\newtheorem{lem}[thm]{Lemma}
\newtheorem{prop}[thm]{Proposition}

\newcommand\dem{\textbf{Proof: }}
\newcommand\Def{\textbf{Definition: }}

\newcommand\rem{\textbf{Remark: }}

\title{Almost reducibility for finitely differentiable $SL(2,\mathbb{R})$-valued quasi-periodic cocycles}
\author{Claire Chavaudret\\
C.R.M. Ennio de Giorgi \\
Collegio Puteano, Scuola Normale Superiore\\
Piazza dei Cavalieri, 3\\
56100 Pisa, Italy \\
e-mail: claire.chavaudret@sns.it}
\date{}

\begin{document}

\maketitle

\textbf{Abstract:} {Quasi-periodic cocycles with a diophantine frequency and with values in SL(2,R) are shown to be almost reducible as long as they are close enough to a constant, in the topology of k times differentiable functions, with k great enough. Almost reducibility is obtained by analytic approximation after a loss of differentiability which only depends on the frequency and on the constant part. As in the analytic case, if their fibered rotation number is diophantine or rational with respect to the frequency, such cocycles are in fact reducible. This extends Eliasson's theorem on Schrödinger cocycles to the differentiable case. }

\section{Introduction}


We are considering quasi-periodic cocycles and the problem of their reducibility in the so-called perturbative case, with emphasis on the two-dimensional case. We mean by quasi-periodic cocycle the fundamental solution of a linear system with quasi-periodic coefficients:

\begin{equation}\label{eq}\forall (t,\theta)\in \mathbb{R}\times \mathbb{T}^d,\ \frac{d}{dt}X^t(\theta)=A(\theta+t\omega)X^t(\theta);\ X^0(\theta)=Id
\end{equation}

\noindent
where $A$ is continuous on the $d$-dimensional torus $\mathbb{T}^d$, matrix-valued and $\omega\in \mathbb{R}^d$ is a rationally independent vector. In this case we say that $X$ is the \textit{cocycle associated to $A$}. 
In this paper we will have a particular interest in the case when $A$ is $sl(2,\mathbb{R})$-valued since in this case it is possible to compute the fibered rotation number of the cocycle and have information on the rotational behaviour of the solutions of \eqref{eq}. 

\bigskip
\noindent It is interesting to define an equivalence relation on cocycles as follow: if $A,B\in C^0(\mathbb{T}^d,gl(n,\mathbb{C}))$, one says that $A$ and $B$ are \textit{conjugated in the sense of cocycles}, or just \textit{conjugated}, if there exists a map $Z$ which is continuous on the torus $2\mathbb{T}^d=\mathbb{R}^d/ 2\mathbb{Z}^d$ such that

\begin{equation}\label{conj}\forall \theta\in 2\mathbb{T}^d, \ \frac{d}{dt}Z(\theta+t\omega)_{\mid t=0}= A(\theta)Z(\theta)-Z(\theta)B(\theta)
\end{equation}

\noindent This kind of conjugation preserves some important dynamical invariants, such as 
the Lyapunov exponents or, if $n=2$, the rotation number which is defined up to $\pi \langle \mathbb{Z}^d,\omega\rangle$. A natural question arises when dealing with a cocycle: can it be conjugated, in the sense of cocycles, to the solution of a system with constant coefficients? When it is so, one says that the solution cocycle is \textit{reducible}. More precisely, a cocycle $X$ solution of \eqref{eq} is reducible if \eqref{conj} holds for some constant $B$. 
In this case we also say that $A$ is \textit{reducible to $B$ by $Z$}, which is equivalent to

\begin{equation}\forall (t,\theta)\in \mathbb{R}\times 2\mathbb{T}^d, \ X^t(\theta)=Z(\theta+t\omega)^{-1}e^{tB}Z(\theta)
\end{equation}

\noindent
Reducibility is an interesting dynamical property since it gives access to most of the information on a cocycle's dynamics: the constant part $B$ will reveal the Lyapunov exponents (real part of the spectrum of $B$) as well as the invariant subbundles of the cocycle and the rotational properties of the system's solutions.

\bigskip
\noindent
It is well known since Floquet's theory that every periodic cocycle (i.e, in the notation above, when $d=1$) is reducible (notice that we have allowed one period doubling in our definition of reducibility). However, the presence of at least two incommensurable frequencies in the coefficient of the system gives rise to non-reducible cocycles. 


\noindent
To go round this difficulty, some authors have introduced a notion of almost reducibility for quasi-periodic cocycles. In some topology $\mathcal{C}$, a cocycle is said to be almost reducible if it can be conjugated, in the sense above, with $Z$ of class $\mathcal{C}$, to another cocycle which is $\mathcal{C}$-arbitrarily close to a constant cocycle.


\noindent
Many results about reducibility and almost reducibility of quasi-periodic cocycles have been obtained in the perturbative case, i.e the case when the vector $\omega$ satisfies a diophantine condition and \eqref{eq} has the following form:

\begin{equation}\label{per1}\frac{d}{dt}X^t(\theta)=(A+F(\theta+t\omega)) X^t(\theta); \ X^0(\theta)=Id
\end{equation}

\noindent
where the coefficient $A+F(\theta+t\omega)$ is close enough to a constant, with a closeness condition related to the diophantine condition on $\omega$: 

\begin{equation}\label{per2}\mid \mid F\mid \mid _\mathcal{C} \leq \epsilon ( n,d,\omega, A)
\end{equation}

\noindent
Then, if $\mathcal{C}$ stands for some analytic class, it is known that
\begin{itemize}
\item every cocycle is almost reducible (\cite{El01}, \cite{C2})
\item almost all cocycles are reducible, when considering a generic one-parameter family (\cite{El92}, \cite{HY} completed with \cite{C1})
\item reducible cocycles are dense (\cite{C2})
\item in the $SO(3)$-valued case, also non reducible cocycles are dense (\cite{El02}). 
\end{itemize}

\noindent
In fact, \cite{El92} also investigates the link between reducibility and the rotational behaviour of the solutions, showing that Schrödinger cocycles are reducible if and only if their fibered rotation number either satisfies a diophantine condition or is rational with respect to $\omega$; this result was extended to general $SL(2,\mathbb{R})$-valued cocycles in \cite{SH}. 

\bigskip
\noindent
Here we shall consider the perturbative case, but in a finitely differentiable topology, a case in which little is known. The aim of this work is to show that in the perturbative regime described by \eqref{per1} and \eqref{per2}, for cocycles which are sufficiently smooth but finitely differentiable, say $C^k$ for some $k\geq k_0(d,\omega,A)$, and have values in $SL(2,\mathbb{R})$, 
all cocycles are almost reducible in a finitely differentiable topology $C^{k-D}$ with a loss of differentiability $D$ which is independent of the initial regularity $k$; in fact, we state this theorem in such a way that it also holds for cocycles with values in other Lie groups. More precisely, we will prove, for $G$ amongst $SL(2,\mathbb{C}),SL(2,\mathbb{R}),O(2), GL(n,\mathbb{C}), U(n)$, letting $\mathcal{G}$ be its Lie algebra:

\begin{thm}\label{PRCkintro} Let $A\in \mathcal{G}$. There exists $k_0,D\in \mathbb{N}$ such that if $k\geq k_0$, there exists 
$\epsilon_0>0$ such that if $F\in C^k(\mathbb{T}^d,\mathcal{G})$ and $\mid \mid F\mid \mid _k\leq \epsilon_0$, then there exist $Z_\infty\in C^{k-D}(\mathbb{T}^d,G)$ and $\bar{A}_\infty\in C^{k-D}(\mathbb{T}^d,\mathcal{G})$ such that

\begin{equation}\forall \theta\in \mathbb{T}^d, \ \partial_\omega Z_\infty(\theta)= (A+F(\theta))Z_\infty (\theta)-Z_\infty(\theta) \bar{A}_\infty(\theta)
\end{equation}

\noindent
and $\bar{A}_\infty$ is the limit, in $C^{k-D}(\mathbb{T}^d,\mathcal{G})$, of reducible functions. 

\end{thm}

\noindent Theorem \ref{PRCkintro} is about almost reducibility of differentiable cocycles. It easily implies density of reducible cocycles near a constant. The reason why it holds 
for those Lie groups is that it is based on another theorem which holds for many classical Lie groups (see \cite{C2}), but we can apply it here only when no period doubling is needed, that is to say, in the complex case or in the 2-dimensional case. The case of symplectic cocycles in dimension greater than 2, for instance, is still an open problem. A similar result was shown in the appendix of \cite{FK} in the case of a $C^\infty$-cocycle with two frequencies ($d=2$) and with values in $SL(2,\mathbb{R})$.

\bigskip
Focusing now on the 2-dimensional case, we will show that given a $SL(2,\mathbb{R})$-valued cocycle, if its fibered rotation number satisfies a diophantine condition or is rational with respect to $\omega$, then it is in fact reducible, thus extending Eliasson's theorem in \cite{El92} to the differentiable case:

\begin{thm}\label{reddiophintro}Let $A\in sl(2,\mathbb{R})$. There exists $k_0,D\in \mathbb{N}$ such that if $k\geq k_0$ and $F\in C^k(\mathbb{T}^d,sl(2,\mathbb{R}))$, there exists 
$\epsilon_0>0$ such that if $\mid \mid F\mid \mid _k\leq \epsilon_0$ and the fibered rotation number $\rho(A+F)$ has the form $2\pi \langle m,\omega\rangle, m\in \mathbb{Z}^d$ or satisfies a diophantine condition with respect to $\omega$:

$$\exists \kappa>0, \ \forall m\in \mathbb{Z}^d\setminus\{0\},\ \mid \rho(A+F)-2\pi \langle m,\omega \rangle \mid \geq \frac{\kappa}{\mid m\mid ^\tau}$$

\noindent
where $\tau$ is the diophantine exponent of $\omega$, then the cocycle associated to $A+F$ is reducible in $C^{k-D}(\mathbb{T}^d,SL(2,\mathbb{R}))$. 

\end{thm}

\bigskip
\noindent
The proof of Theorem \ref{PRCkintro} relies essentially on a proposition shown in \cite{C2}, which was used as an inductive lemma in a KAM scheme to show almost reducibility for some analytic and Gevrey cocycles. Here, we use it to get a good control on a sequence of analytic cocycles which, following an idea by Zehnder (\cite{Z}), are constructed in such a way that they approach a given differentiable cocycle. Since they are shown to be conjugated to something which becomes closer and closer to a constant, one finally gets almost reducibility for their limit in a topology with a finite loss of differentiability with respect to the initial topology. 

\bigskip
\noindent
The specificity of $SL(2,\mathbb{R})$, however, is that the eigenvalues of the constant part get closer to 0 every time that, in the KAM scheme, a resonance is removed. Thus, non reducibility implies that the fibered rotation number of the limit cocycle cannot be diophantine, and so, by invariance through conjugation in the sense of cocycles, neither can the fibered rotation number of the initial cocycle, which gives Theorem \ref{reddiophintro}. We then easily get an application to Schrödinger cocycles inspired by \cite{El92}.

\subsection{Definitions and assumptions}

\noindent
Throughout this paper, we will make the following assumption.

\bigskip
\textbf{Assumption:} There exist $0<\kappa<1$ and $ \tau \geq \max (1,d-1)$ such that

\begin{equation}\forall m\in \mathbb{Z}^d\setminus\{0\}, \mid \langle m,\omega \rangle \mid \geq \frac{\kappa}{\mid m\mid ^\tau}.
\end{equation}

\noindent The numbers $\kappa$ and $\tau$ will be fixed from now on. 
This is a diophantine condition on $\omega$. We shall define other types of diophantine conditions, which refer to the vector $\omega$.

\bigskip
\Def Let $z\in \mathbb{R}$; we say that $z$ is diophantine with respect to $\omega\in \mathbb{R}^d$ and we write $z\in DC_\omega$ if there exists $\kappa'>0,\tau'>\max (1,d-1)$ such that for all $m\in \mathbb{Z}^d\setminus\{0\}$,

\begin{equation}\mid z-2\pi \langle m,\omega\rangle\mid \geq \frac{\kappa'}{\mid m\mid ^{\tau'}}.
\end{equation}

\noindent We will also denote by $DC_\omega(\tau)$ the set of numbers $z\in \mathbb{R}$ such that there exists $\kappa'>0$ satisfying, for all $m\in \mathbb{Z}^d\setminus\{0\}$,

\begin{equation}\mid z-2\pi \langle m,\omega\rangle\mid \geq \frac{\kappa'}{\mid m\mid ^{\tau}}.
\end{equation}

\noindent
The following diophantine condition is also known as "second Melnikov condition" and refers to the spectrum of a matrix:

\bigskip
\Def Let $A\in gl(n,\mathbb{C})$ and $\{ \alpha_1,\dots ,\alpha_n\}$ its spectrum. Let $\kappa'>0, N\in \mathbb{N}$; we say that $A$ \textit{has a $DC^N_\omega (\kappa',\tau)$ spectrum} if 

\begin{equation}\forall 1\leq j,k\leq n,\forall m\in \mathbb{Z}^d\setminus\{0\}, \mid m\mid \leq N\Rightarrow \mid \operatorname{Im}(\alpha_j)-\operatorname{Im}(\alpha_k) -2\pi \langle m,\omega \rangle \mid \geq \frac{\kappa'}{\mid m\mid ^\tau}.
\end{equation}

\noindent 
If $A\in sl(2,\mathbb{R})$ with spectrum $\{\pm \alpha\}$, this reduces to

\begin{equation}\forall m\in \mathbb{Z}^d\setminus\{0\}, \mid m\mid \leq N\Rightarrow \mid 2\operatorname{Im}(\alpha) -2\pi \langle m,\omega \rangle \mid \geq \frac{\kappa'}{\mid m\mid ^\tau}.
\end{equation}

\Def We will denote by $\mathcal{M}_\omega$ the set of numbers which are rational with respect to $\omega$, i.e

\begin{equation}\label{momega}\mathcal{M}_\omega=\{ 2\pi \langle m,\omega\rangle,\ m\in \mathbb{Z}^d\}.
\end{equation}

\noindent It has a module structure, therefore it is sometimes called the \textit{frequency module}.

\bigskip
\noindent 
Now we recall the definition of the fibered rotation number of an $SL(2,\mathbb{R})$-valued cocycle: 



\bigskip
\Def
 Let $A\in C^k(\mathbb{T}^d,sl(2,\mathbb{R}))$. We will denote by $\rho (A)$ and refer to as the \textit{fibered rotation number of the cocycle $X$ associated to $A$} the number

$$\rho(A)= \lim_{t \rightarrow +\infty} \frac{1}{t} Arg (X^t(\theta)\phi)$$

\noindent 
where $Arg$ is the variation of the argument. 


\bigskip
\rem \begin{itemize}
\item This number does not depend on the choice of $\phi$ (see \cite{El92}, appendix);
\item If $A$ and $B$ are conjugated in the sense of cocycles, then $\rho(A)=\rho(B)+\langle m,\omega\rangle$ for some $m\in \frac{1}{2}\mathbb{Z}^d$;
\item if $A$ is reducible to $B$ by some function $Z$, then $\rho(A)$ coincides with the Floquet exponent of $A$ i.e the modulus of the imaginary part of the eigenvalues of $B$ (which is well-defined only modulo $\frac{1}{2}\mathcal{M}_\omega$).
\end{itemize}

\bigskip
\noindent
\textbf{Notations:} 
The usual operator norm will be denoted by $\mid \mid .\mid \mid$. In the space $C^k(2\mathbb{T}^d,gl(n,\mathbb{C}))$ of $k$ times differentiable matrix-valued functions on the double torus, we will 
use the norm

$$\mid \mid F\mid \mid _k= \sup_{k'\leq k; \theta\in \mathbb{T}^d} \mid \mid d^{k'} F (\theta)\mid \mid .$$

For any map $Z\in C^1(2\mathbb{T}^d,gl(n,\mathbb{C}))$ we will denote by $\partial_\omega Z$ the derivative of $Z$ in the direction $\omega$: 

$$\forall \theta \in 2\mathbb{T}^d, \partial_\omega Z(\theta)=\frac{d}{dt}Z(\theta+t\omega)_{\mid t=0}.$$

\section{A lemma on analytic cocycles}

We first recall a proposition which will be used as inductive step in the proof of Theorem \ref{PRCkintro}. It was proved in \cite{C2} (Proposition 2.14). 

\bigskip
\noindent \textbf{Notations:}
In the following, for $r>0$ and any set $E$, we will denote by $C^\omega_r(\mathbb{T}^d,E)$ the space of functions which are analytic on a "strip" 

$$\{z=(z_1,\dots,z_d)\in \mathbb{C}^d, \mid \operatorname{Im} z_1
\mid <r,\dots , \mid {Im} z_d\mid <r\}$$

\noindent and are 1-periodic in $\operatorname{Re} z_1,\dots , \operatorname{Re} z_d$ and whose restriction on $\mathbb{R}^d$ has values in $E$. The writing $C^\omega_r(2\mathbb{T}^d,E)$ will stand for functions which are analytic on a strip and $E$-valued on $\mathbb{R}^d$, but only 2-periodic in $\operatorname{Re} z_1,\dots , \operatorname{Re} z_d$.

\bigskip
\noindent The norm in $C^\omega_r(2\mathbb{T}^d,gl(n,\mathbb{C}))$ will be written $\mid .\mid _r$.

\bigskip
\noindent We shall fix a Lie group $G$ amongst $GL(n,\mathbb{C}), U(n), SL(2,\mathbb{C}),SL(2,\mathbb{R}),O(2)$ and denote by $\mathcal{G}$ its Lie algebra. 

\bigskip
\noindent 
To simplify the statements, we shall use the following technical abreviations:

\begin{equation}\label{kappa''}\left\{\begin{array}{c}
N(r,\epsilon)=\frac{1}{2\pi r}|\log \epsilon|\\
R(r,r')= \frac{1}{(r-r')^8}80^{4}(\frac{1}{2}n(n-1)+1)^2\\
\kappa''(r,r',\epsilon)=\frac{\kappa}{n(8R(r,r')^{\frac{1}{2}n(n-1)+1}{N(r,\epsilon)})^{\tau}}
\end{array}\right.\end{equation}
%

%


%





\begin{prop}\label{iter3}
Let
\begin{itemize}
\item $A\in \mathcal{G}$,
\item $r\leq \frac{1}{2},r''\in [\frac{95}{96}r,r[$, 
\item $\bar{A},\bar{F}\in C^\omega_{r}(\mathbb{T}^d,\mathcal{G})$ and $\Psi\in C^\omega_{r}(2\mathbb{T}^d,G)$, 
\item $\epsilon=|\bar{F}|_{r}$,
\end{itemize}

\noindent There exists $C>0$ depending only on $n,d,\kappa,\tau$ and there exists $D\in\mathbb{N}$ depending only on 
$n,d,\tau$ 
such that if

\begin{enumerate}

\item $\bar{A}$ is reducible to $A$ 
by $\Psi$, 

\item $\Psi$ is such that for all $H$ continuous on $\mathbb{T}^d$, $\Psi H\Psi^{-1}$ is continuous on $\mathbb{T}^d$,

\item

\begin{equation}\label{epsA+} \epsilon 
\leq \frac{C}{(||A||+1)^{D}}(r-r'')^{D},
\end{equation}

\item $|\Psi|_{r}\leq (\frac{1}{\epsilon})^{-\frac{1}{2}(r-r'')}$ and
$|\Psi^{-1}|_{r}\leq (\frac{1}{\epsilon})
^{-\frac{1}{2}(r-r'')}$, 


\end{enumerate}
 
\noindent then there exist
\begin{itemize}
\item $\epsilon'\in [\epsilon^{R(r,r'')^{n^2}}, \epsilon^{100}]$;
\item $Z'\in C^\omega_{r''}(\mathbb{T}^d,G)$, 
\item $\bar{A}',\bar{F}'\in C^\omega_{r''}(\mathbb{T}^d,\mathcal{G})$, 
\item $\Psi'\in C^\omega_{r}(2\mathbb{T}^d, G)$, 
\item ${A}'\in \mathcal{G}$ 

\end{itemize}
\noindent satisfying the following properties:

\begin{enumerate}

\item \label{5} $\bar{A}'$ is reducible by $\Psi'$ to ${A}'$, 
\item $\Psi'$ is such that for all $H$ continuous on $\mathbb{T}^d$, $\Psi' H\Psi^{'-1}$ is continuous on $\mathbb{T}^d$,

\item \label{8+} $|\bar{F}'|_{r''} \leq  \epsilon '$,

\item \label{ren+} $|\Psi'|_{r''}\leq (\frac{1}{\epsilon'})^{\frac{1}{4}(r-r'')}$ and
$|\Psi'^{-1}|_{r''}\leq 
(\frac{1}{\epsilon'})^{\frac{1}{4}(r-r'')}$,

\item \label{0} $||A'||\leq ||A||+\mid \log \epsilon\mid\left(\frac{1}{r-r''}\right) ^{D}$;

\item \label{4+}

\begin{equation}\partial_\omega Z'=(\bar{A}+\bar{F})Z'
-Z'(\bar{A}'+\bar{F}'),\end{equation}

\item \label{6} 

\begin{equation}|Z'-Id|_{r''}\leq \frac{1}{C}\left(\frac{(1+||A||)|\log\epsilon|}
{r-r''}\right)^{D}
\epsilon^{1-4(r-r'')}
\end{equation}

\noindent and so does $(Z')^{-1}-Id$.
%


\end{enumerate}

\noindent Moreover, 
\begin{itemize}
\item if $\mathcal{G}=o(2)$ or $u(n)$, the same holds with the weaker condition

\begin{equation}\label{epsorth+} \epsilon 
\leq C(r-r'')^{D}
\end{equation}

\noindent instead of \eqref{epsA+};

\item if $\mathcal{G}=sl(2,\mathbb{C})$ or $sl(2,\mathbb{R})$, then either $\Psi'^{-1}\Psi$ is the identity or $\mid \mid A'\mid \mid\leq \kappa''(r,r'',\epsilon) +\epsilon^{\frac{1}{2}}$.

\end{itemize}

\end{prop}

\bigskip
The proof of Proposition \ref{iter3} given in \cite{C2} also implies the following:

\bigskip
\noindent 
If $A$ has a $DC^N_\omega(\kappa'',\tau)$ spectrum with $N=N(r,\epsilon) $ and $\kappa''=\kappa''(r,r'',\epsilon)$, then $\Psi=\Psi'$.

\bigskip
By construction, functions $\Psi$ and $\Psi'$ also satisfy the following, in case $\mathcal{G}=sl(2,\mathbb{R})$:

\bigskip
\textit{If $\Psi$ satisfies}: 

\bigskip
for all $A,A'\in C^0(\mathbb{T}^d,sl(2,\mathbb{R}))$, $\partial_\omega \Psi=A\Psi-\Psi A'\Rightarrow \rho(A)=\rho(A')+2\pi \langle m,\omega\rangle$ for some $m\in \frac{1}{2} \mathbb{Z}^d$ 

\bigskip
\textit{then $\Psi'$ satisfies the same property:}

\bigskip
for all $A,A'\in C^0(\mathbb{T}^d,sl(2,\mathbb{R}))$, $\partial_\omega \Psi'=A\Psi'-\Psi' A'\Rightarrow \rho(A)=\rho(A')+2\pi \langle m',\omega\rangle$ for some $m'\in \frac{1}{2}\mathbb{Z}^d$ 

\bigskip
\textit{with $m=m'$ if and only if $\Psi=\Psi'$. }

\bigskip
\noindent Moreover, looking more closely at the proof of Proposition 2.1 in \cite{C2}, one has for $\mathcal{G}=sl(2,\mathbb{C})$ or $sl(2,\mathbb{R})$:

\begin{equation}\label{nbrot}\mid \mid \Psi'^{-1}\Psi (A \Psi^{-1}\Psi'- \partial_\omega( \Psi^{-1}\Psi')  )-A' \mid  \mid \leq \sqrt{\epsilon}
\end{equation}

\noindent therefore (since the rotation number of a constant cocycle is given by the imaginary part of its spectrum)

\begin{equation}\label{nbrot'}\mid \rho(A)-\rho(A')+2\pi \langle m'-m,\omega\rangle \mid \leq \sqrt{\epsilon}.
\end{equation}

\section{Almost reducibility}

First we need a numerical lemma:

\begin{lem}\label{num}Let $C>0,D,k\in \mathbb{N}$,$\epsilon_j=\frac{1}{j^k}$. There exists $k_1\in\mathbb{N}\setminus\{0\}$ such that if $k\geq k_1$, then for all $j\geq 2$,

\begin{equation}\label{borne}C[j(j+1)\mid \log \epsilon_j\mid ]^D \epsilon_j^{1-\frac{4}{j(j+1)}}\leq \frac{1}{(j+1)^2}.
\end{equation}

\end{lem}

\dem Equation \eqref{borne} is equivalent to

\begin{equation}\label{borne'}C[j(j+1)k\log j ]^D (\frac{1}{j})^{k(1-\frac{4}{j(j+1)})}\leq \frac{1}{(j+1)^2}.
\end{equation}

There exists $k_1$ such that for all $k\geq k_1,j\geq 2$, 

\begin{equation}C\frac{k^D}{j^{(\frac{k}{2}-3D)  }}\leq \frac{1}{(j+1)^2} 
\end{equation}

so \eqref{borne} holds. $\Box$

\bigskip
\noindent We will now state the main result for $G$ among $GL(n,\mathbb{C}),U(n), SL(2,\mathbb{C}),SL(2,\mathbb{R}),O(2)$. We shall denote by $\mathcal{G}$ the Lie algebra associated to $G$. 

\begin{thm}\label{PRck} Let $A\in \mathcal{G}$. There exists
$k_0,D'\in \mathbb{N}$ 
only depending on $n,d,\kappa,\tau,A$
such that for all $k\geq k_0$ and 
$F\in C^k(\mathbb{T}^d,\mathcal{G})$, there exists $\epsilon_0$ depending only on $n,d,\kappa,\tau,A,k$ such that if $||F||_{k}\leq \epsilon_0$, then 
there exist 
\begin{itemize}
\item $Z_\infty\in C^{k-D'}(\mathbb{T}^d,G)$, 
\item $\bar{A}_\infty\in C^{k-D'}(\mathbb{T}^d,\mathcal{G})$,
\item a sequence $(\bar{A}_j)_{j\geq 1}$ of functions in $C^{k-D'}(\mathbb{T}^d,\mathcal{G})$,
\item a sequence $(\Psi_j)_{j\geq1}$ of functions in $C^{k-D'}(2\mathbb{T}^d,G)$,
\item a sequence $(A_j)_{j\geq 1}$ of elements of $\mathcal{G}$
\end{itemize}
such that 
\begin{enumerate}
\item \label{limite} $\bar{A}_\infty$ is the limit in $C^{k-D'}
(\mathbb{T}^d,\mathcal{G})$ 
of the sequence $\bar{A}_j$, 
\item \label{rÈductibles} the functions $\bar{A}_j$ are reducible to $A_j$ by $\Psi_j$, 

\item \label{conjuguÈs}

\begin{equation}\partial_\omega Z_\infty(\theta)=(A+F(\theta))Z_\infty(\theta)
-Z_\infty(\theta)\bar{A}_\infty(\theta).
\end{equation}

\end{enumerate}

%

Moreover, in the case $\mathcal{G}=sl(2,\mathbb{R})$, 
there exist
\begin{itemize}
\item a sequence $(M_j)_{j\geq 1}$ of elements of $\frac{1}{2}\mathbb{Z}^d$ 
\item unbounded sequences $(N_j)_{j\geq 1}$ and $(R_j)_{j\geq 1}$ of integers
\end{itemize}
such that letting $\kappa_j=\frac{\kappa}{2(8R_j^2N_j)^\tau}$, then
\begin{itemize}

\item the sequence $(\kappa_j)$ is summable,


\item for all $A_1,A_2\in C^0(\mathbb{T}^d,sl(2,\mathbb{R}))$ and all $j$, 

$$\partial_\omega \Psi_j=A_1\Psi_j-\Psi_jA_2\Rightarrow \rho(A_1)=\rho(A_2)+2\pi \sum_{l=1}^j \langle M_l,\omega\rangle,$$

\item either $A_j$ has a $DC^{N_j}_\omega (\kappa_j,\tau)$ spectrum,
 which implies that $M_j=0$ or $M_{j-1}=0$, or $R_{j-1}N_{j-1}<\mid M_j\mid \leq N_j$ and in that case $\sigma(A_j)\subset B(0,\kappa_{j-1})$;
 
 \item $\mid \rho(A_{j+1})-(\rho(A_j)+2\pi \langle M_j,\omega\rangle)\mid \leq \kappa_j$;
 
 \item if there exists $J\geq 1$ such that $M_j=0$ for all $j\geq J$, then $A+F$ is reducible.
 \end{itemize}
\end{thm}

\bigskip
\noindent In this statement, properties \ref{limite}, \ref{rÈductibles} and \ref{conjuguÈs} are sufficient to get Theorem \ref{PRCkintro}, but the other properties will be used in the application to $SL(2,\mathbb{R})$-valued cocycles. 

\bigskip
\dem $\bullet$ By \cite{Z}, there exists a sequence $(F_j)_{j\geq 1}$, 
$F_j\in C^\omega_\frac{1}{j}(\mathbb{T}^d,\mathcal{G})$ and a universal constant $C'$, such that

\begin{equation}\label{suite}\left\{\begin{array}{c}
||F_j-F||_k\rightarrow 0\ \mathrm{when}\ j\rightarrow +\infty\\
|F_j|_\frac{1}{j}\leq C'||F||_k\\
|F_{j+1}-F_j|_\frac{1}{j+1}\leq C'(\frac{1}{j})^k
||F||_k\\
\end{array}\right.\end{equation}
%



 
\bigskip
\noindent Moreover, this sequence is obtained from $F$ regardless of its regularity, i.e if 
$k\leq k'$ and $F\in C^{k'}$, then properties \eqref{suite} hold with $k'$ instead of $k$ 
(since $F_j$ is the convolution of $F$ with a map which does not depend on $k$).

\bigskip
\noindent Let $C>0,D\in \mathbb{N}$ be as in Proposition \ref{iter3}. One can assume $C\leq \frac{1}{2}$. Recall that these numbers only depend on $n,d,\kappa,\tau,A$. 
For all $r>r'>0$, let


$$\epsilon'_0(r,r')=C(r-r')^{D}.$$ 



\noindent
Let $k_1$ be as in Lemma \ref{num} and let $k_0\geq k_1$ be a number depending only on $n,d,\kappa,\tau,A$ such that for all $j\geq 2$,

$$\frac{C}{j^{k_0}}\leq \epsilon'_0(\frac{1}{j},\frac{1}{j+1}).$$

\noindent
Assume $k\geq k_0$ and let for all $j\geq 1$,

$$\epsilon'_j=\frac{C}{j^k}$$


\noindent so that $\epsilon_j'\leq\epsilon'_0(\frac{1}{j},\frac{1}{j+1})$, and let

$$\alpha_j=\frac{4}{j(j+1)}.$$
 





\bigskip
$\bullet$ \textit{First step:} Assume that

\begin{equation}\label{fk}C'||F||_k
\leq \epsilon'_2=\frac{C}{2^k}\end{equation}

\noindent (notice that this condition on $||F||_k$ only depends on $n,d,\kappa,\tau,A,k$). Then

$$|F_2|_\frac{1}{2}
\leq \epsilon'_2\leq \epsilon'_0(\frac{1}{2},\frac{1}{3})$$
 
\noindent therefore, by Proposition \ref{iter3}, there exist
\begin{itemize}
\item 
$\epsilon''_2\leq |F_2|_\frac{1}{2}^{100}$
\item $Z_2\in C^\omega_\frac{1}{3}(\mathbb{T}^d,
G)$, 
\item $\bar{A}_2\in C^\omega_\frac{1}{3}(\mathbb{T}^d,
\mathcal{G})$ 
\item $\bar{F}_2\in C^\omega_\frac{1}{3}(\mathbb{T}^d,
\mathcal{G})$ 
 
\item $\Psi_2\in C^\omega_1(2\mathbb{T}^d,G)$
\item 
$A_2\in \mathcal{G}$
\end{itemize}
such that

\begin{enumerate}

\item

$\bar{A}_2$ is reducible to $A_2$ by $\Psi_2$,

\item 
$\mid \bar{F}_2\mid_\frac{1}{3}\leq\epsilon''_2\leq \frac{1}{2}\epsilon'_3$,

\item $\mid\Psi_2\mid_\frac{1}{3}\leq \left(\frac{1}{\epsilon''_2}\right)^{\alpha_2}$, as well as $\Psi_2^{-1}$,


\item 
\begin{equation} \partial_\omega Z_2=(A+F_2)Z_2-Z_2(\bar{A}_2
+\bar{F}_2),
\end{equation}

\item

\begin{equation} \label{zm0} |Z_2-Id|_{\frac{1}{3} }\leq 
\frac{1}{C} \left( \frac{4}{\alpha_2}(1+||A||)|\log |F_2|_\frac{1}{2}| \right)^{D}
 |F_2|_\frac{1}{2}^{1-\alpha_2}
\end{equation}

\noindent and $Z_2^{-1}$ satisfies the same estimate,

\item and if $\mathcal{G}=sl(2,\mathbb{R})$, $\Psi_2$ satisfies: for all $A,A'\in C^0(\mathbb{T}^d,sl(2,\mathbb{R}))$, 

\begin{equation}\partial_\omega \Psi_2 = A\Psi_2 -\Psi_2 A'\Rightarrow \exists M_1\in \frac{1}{2}\mathbb{Z}^d, \rho(A)=\rho(A')+2\pi \langle M_1,\omega\rangle.
\end{equation}

\end{enumerate}

\noindent Property \eqref{zm0} implies

\begin{equation} \label{zm0'} \begin{split}|Z_2-Id|_{\frac{1}{3} }&\leq 
\frac{1}{C} \left( \frac{4}{\alpha_2}(1+||A||)|\log  \epsilon'_2| \right)^{D}
 ( \epsilon'_2)^{1-\alpha_2}\\
\end{split}\end{equation}


%

\noindent and $Z_2^{-1}$ satisfies the same estimate as $Z_2$. Lemma \ref{num} then implies that $|Z_2-Id|_{\frac{1}{3} }\leq \frac{1}{9}$.

\bigskip
\noindent $\bullet$ \textit{Induction step:} Let $j\geq 2$. Suppose that there exists
\begin{itemize}
\item $A_j\in \mathcal{G}$ 
\item $\Psi_j\in C^\omega_1(2\mathbb{T}^d,G)$
\item 
$\bar{A}_j\in C^\omega_\frac{1}{j+1}(\mathbb{T}^d,
\mathcal{G})$ and $\bar{F}_j\in C^\omega_\frac{1}{j+1}(\mathbb{T}^d,
\mathcal{G})$
\item $\bar{Z}_j\in C^\omega_\frac{1}{j+1}(\mathbb{T}^d,
G)$, 
\end{itemize}
 such that
 
 \begin{enumerate}
 
 \item $\bar{A}_j$ is reducible to $A_j$ by $\Psi_j$,
 
 \item $\mid \Psi_j\mid_\frac{1}{j+1} \leq |\bar{F}_j|_\frac{1}{j+1}^{-\alpha_{j}}$,
 
 \item $|\bar{F}_j|_\frac{1}{j+1}\leq \frac{1}{2} \epsilon'_{j+1}$,

\item 
$
|\bar{Z}_j-Id|_\frac{1}{j+1}\leq  2\sum_{l=3}^{j+1}\frac{1}{l^2}\leq 1$,

\item if $\mathcal{G}=sl(2,\mathbb{R})$, for all $A,A'\in C^0(\mathbb{T}^d,sl(2,\mathbb{R}))$, 

\begin{equation}\partial_\omega \Psi_j=A\Psi_j -\Psi_j A' \Rightarrow \exists M_1,\dots, M_{j-1}\in \frac{1}{2}\mathbb{Z}^d, \rho(A')=\rho(A)+2\pi \langle M_1+\dots +M_{j-1},
\omega\rangle 
\end{equation}

\item and

\begin{equation}\label{zjbar} \partial_\omega \bar{Z}_j=(A+F_j)\bar{Z}_j-\bar{Z}_j(\bar{A}_j
+\bar{F}_j).
\end{equation}

\end{enumerate}

\noindent Then 

\begin{equation} \partial_\omega \bar{Z}_j=(A+F_{j+1})\bar{Z}_j-\bar{Z}_j(\bar{A}_j
+\bar{Z}_j^{-1}(F_{j+1}-F_{j})\bar{Z}_j+\bar{F}_j)
\end{equation}

\noindent and moreover, 
by \eqref{suite}, 

\begin{equation}|\bar{Z}_j^{-1}(F_{j+1}-F_j)\bar{Z}_j+\bar{F}_j|_\frac{1}{j+1}\leq \frac{1}{2}\epsilon'_{j+1}
+\frac{4}{j^k}C'||F||_k\end{equation}

\noindent which implies, by assumption \eqref{fk}, that

\begin{equation}|\bar{Z}_j^{-1}(F_{j+1}-F_j)\bar{Z}_j+\bar{F}_j|_\frac{1}{j+1}
\leq \epsilon'_{j+1}\leq \epsilon'_0(\frac{1}{j+1},\frac{1}{j+2})\\
\end{equation}

\noindent so one can apply Proposition \ref{iter3}: denoting $\tilde{\epsilon}_j=|\bar{Z}_j^{-1}(F_{j+1}-F_{j})\bar{Z}_j+\bar{F}_j|_\frac{1}{j+1}$, there exists 
\begin{itemize}
\item $\epsilon''_{j+1}\leq \tilde{\epsilon}_j^{100}$
\item 
$Z_{j+1}\in C^\omega_\frac{1}{j+2}(\mathbb{T}^d,
G)$, 
\item 
$\bar{A}_{j+1}\in C^\omega_\frac{1}{j+2}(\mathbb{T}^d,
\mathcal{G})$ and $\bar{F}_{j+1}\in C^\omega_\frac{1}{j+2}
(\mathbb{T}^d,
\mathcal{G})$ 
\item $\Psi_{j+1}\in C^\omega_1(2\mathbb{T}^d,G)$
\item $A_{j+1}\in \mathcal{G}$ 
\end{itemize}
such that 

\begin{enumerate}

\item $\bar{A}_{j+1}$ is reducible by $\Psi_{j+1}$ to $A_{j+1}$,

\item 
$|\bar{F}_{j+1}|_\frac{1}{j+2}\leq \epsilon''_{j+1}\leq (\epsilon'_{j+1})^{100}\leq \frac{1}{2}\epsilon'_{j+2}$,

\item $\mid \Psi_{j+1}\mid_{\frac{1}{j+2}}\leq (\epsilon_{j+1}'')^{-\alpha_{j+1}}\leq |\bar{F}_{j+1}|_\frac{1}{j+2}^{-\alpha_{j+1}}$,


\item 

\begin{equation} \partial_\omega Z_{j+1}=(\bar{A}_j
+\bar{Z}_j^{-1}(F_{j+1}-F_{j})\bar{Z}_j+\bar{F}_j)Z_{j+1}-Z_{j+1}
(\bar{A}_{j+1}+\bar{F}_{j+1}),
\end{equation}

\item 

\begin{equation}\label{zm} |Z_{j+1}-Id|_\frac{1}{j+2}\leq \frac{1}{C} \left( \frac{4}{\alpha_{j+1}}(1+||A||)|\log \tilde{\epsilon}_j| \right)^{D}
(\tilde{\epsilon}_j)^{1-\alpha_{j+1} },
 \end{equation}

\item if $\mathcal{G}=sl(2,\mathbb{R})$, for all $A,A'\in C^0(\mathbb{T}^d,sl(2,\mathbb{R}))$, 

\begin{equation}\partial_\omega \Psi_{j+1}=A\Psi_{j+1} -\Psi_{j+1} A' \Rightarrow \exists M_1,\dots, M_{j}\in \frac{1}{2}\mathbb{Z}^d, \rho(A')=\rho(A)+2\pi \langle M_1+\dots +M_{j},
\omega\rangle ,
\end{equation}

 \item if $A_j$ has a $DC^{N_j}_\omega(\kappa_j,\tau)$ spectrum, with $N_j=\frac{j+1}{2\pi} \mid \log \tilde{\epsilon}_j\mid $ and $\kappa_j=\frac{\kappa}{2[8R_j^2N_j]^\tau}$ with $R_j=4((j+1)(j+2))^880^4$, then $\Psi_{j+1}=\Psi_j$.

\item $\mid \rho(A_{j+1})-\rho(A_j)+2\pi\langle M_j,\omega\rangle )\mid \leq \kappa_j$.

\end{enumerate}

\noindent Property \eqref{zm} implies

\begin{equation}\label{zm'}\begin{split} |Z_{j+1}-Id|_\frac{1}{j+2} &\leq\frac{(1+\mid \mid A\mid \mid)^D}{C} \left( (j+2)^2|\log \tilde{\epsilon}_j| \right)^{D}
(\tilde{\epsilon}_j)^{1-\alpha_{j+1} }
\end{split} \end{equation}

 
 \noindent so by lemma \ref{num}, $|Z_{j+1}-Id|_\frac{1}{j+2}\leq  \frac{1}{(j+2)^2}$.
 Let $\bar{Z}_{j+1}=\bar{Z}_j Z_{j+1}$. Then

\begin{equation}\label{zmbar} |\bar{Z}_{j+1}-Id|_\frac{1}{j+2}
\leq \mid \bar{Z}_{j}\mid_\frac{1}{j+2}  \mid Z_{j+1}-Id\mid_\frac{1}{j+2} +  \mid \bar{Z}_j-Id\mid_\frac{1}{j+2} 
\leq 2\sum_{l=3}^{j+2}\frac{1}{l^2}.
 \end{equation}

\noindent Property \eqref{zm} also implies

\begin{equation}\label{zm'}\begin{split} |Z_{j+1}-Id|_\frac{1}{j+2} &\leq\tilde{C} \left( (j+2)^2|\log \tilde{\epsilon}_j| \right)^{D}
(\tilde{\epsilon}_j)^{1-\alpha_{j+1} }\\
&\leq \tilde{C} \left( (j+2)^2|\log (\epsilon'_{j+1})| \right)^{D}
(\epsilon'_{j+1})^{1-\alpha_{j+1} }\\
&\leq \tilde{C}' k ^{2D}
(j+1)^{k(\alpha_{j+1} -1)+3D}\\.
\end{split} \end{equation}

%






\bigskip
\noindent $\bullet$ \textit{Conclusion:} So for all $j\geq 2$, there exist
\begin{itemize}
\item
$\bar{Z}_j,Z_j\in C^\omega_\frac{1}{j+1}(\mathbb{T}^d,
G)$, 
\item 
$\bar{A}_j\in C^\omega_\frac{1}{j+1}(\mathbb{T}^d,
\mathcal{G})$ and $\bar{F}_j\in C^\omega_\frac{1}{j+1}(\mathbb{T}^d,
\mathcal{G})$
\item $\Psi_j\in C^\omega_1(2\mathbb{T}^d,G)$,
\item $A_j\in \mathcal{G}$ 
\end{itemize}
such that 

\begin{enumerate}

\item $\bar{Z}_j=Z_1\dots Z_j$,

\item $\bar{A}_j$ is reducible to $A_j$ by $\Psi_j$,

\item $|\bar{F}_j|_\frac{1}{j+1}\leq \epsilon'_{j+1}$,

\item $\mid\Psi_j\mid_{\frac{1}{j+1}}\leq  |\bar{F}_j|_\frac{1}{j+1}^{-\alpha_{j}} $,


\item 

\begin{equation}\label{redck} \partial_\omega \bar{Z}_j(\theta)=(A+F_j(\theta))\bar{Z}_j(\theta)-\bar{Z}_j(\theta)(\bar{A}_j(\theta)
+\bar{F}_j(\theta)),
\end{equation}

\item

\begin{equation} |\bar{Z}_{j}-Id|_\frac{1}{j+1}\leq 1;\   |\bar{Z}_{j}^{-1}-Id|_\frac{1}{j+1}\leq 1
\end{equation}

\item and

\begin{equation}\label{zm''}\begin{split} |Z_{j}-Id|_\frac{1}{j+1} &
\leq \tilde{C}' k^{2D}j^{k(\alpha_{j} -1)+3D}.
\end{split} \end{equation}

\end{enumerate}

Moreover, in the case $G=SL(2,\mathbb{R})$, there exists a sequence $(M_j)_{j\geq 1}$ of elements of $\frac{1}{2}\mathbb{Z}^d$ such that for all $j\geq 1$, 
\begin{itemize}
\item if $A_j$ has a $DC^{N_j}_\omega(\kappa_j,\tau)$ spectrum, with $N_j=\frac{j+1}{2\pi} \mid \log \tilde{\epsilon}_j\mid $ and $\kappa_j=\frac{\kappa}{2[8R_j^2N_j]^\tau}$ with $R_j=4((j+1)(j+2))^880^4$, then $\Psi_{j+1}=\Psi_j$ and $M_j=0$;
\item for all $A,A'\in C^0(\mathbb{T}^d,sl(2,\mathbb{R}))$, 

\begin{equation}\label{mj}\partial_\omega (\Psi_{j+1}\Psi_j^{-1}) = A\Psi_{j+1}\Psi_j^{-1}-\Psi_{j+1}\Psi_j^{-1}A'\Rightarrow \rho(A)=\rho(A')+2\pi\langle M_j,\omega\rangle 
\end{equation}

and either $M_j=0$ (if and only if $\Psi_j=\Psi_{j+1}$), or $M_{j-1}=0$, or

$$R_{j-1}N_{j-1}<M_j\leq N_j;$$

\item $\mid \rho(A_{j+1})-\rho(A_j)+2\pi\langle M_j,\omega\rangle )\mid \leq \kappa_j$.

%

%

\end{itemize}


\bigskip
\noindent $\bullet$ \textit{Convergence:} Now we have to compute the topology in which the sequence $(\bar{Z}_j)$ defined above is Cauchy.
Since 

\begin{equation}\begin{split}|\bar{Z}_j-\bar{Z}_{j+1}|_\frac{1}{j+2}&\leq |\bar{Z}_j|_\frac{1}{j+1}
|Z_{j+1}-Id|_\frac{1}{j+2}\\
&\leq \tilde{C}'k^{2D}
(j+1)^{k(\alpha_{j+1} -1)+3D}\\
\end{split}
\end{equation}

\noindent then for all $k'\in\mathbb{N}$,

\begin{equation}\begin{split}||\bar{Z}_j-\bar{Z}_{j+1}||_{k'}
\leq C_3 (j+1)^{k(\alpha_{j+1} -1)+3D+k'+1}
\end{split}\end{equation}

\noindent for some $C_3$ independent of $j$,
so the sequence $(\bar{Z}_j)$ is Cauchy in the $C^{k'}$ topology if there exists an $j$ such that for all $j'\geq j$,

\begin{equation}\label{cauchy}
k'+1+k(\alpha_{j'+1} -1)+3D<0.
\end{equation}

Let $k'=k-3D-2$. If $j'>4k$, then \eqref{cauchy} holds, therefore $(\bar{Z}_j)$ is Cauchy in the $C^{k'}$ topology. 



\bigskip
\noindent Let $Z_\infty$ be the limit of $(\bar{Z}_j)$ 
in the $C^{k'}$ topology. 
Taking the $C^{k'}$-limit in \eqref{redck}, one gets 

\begin{equation}\partial_\omega Z_\infty(\theta)=(A+F(\theta))Z_\infty(\theta)
-Z_\infty(\theta)\bar{A}_\infty(\theta)
\end{equation}

\noindent where $\bar{A}_\infty\in C^{k'}(\mathbb{T}^d,\mathcal{G})$ is the limit in 
$C^{k'}(\mathbb{T}^d,\mathcal{G})$ of the reducible functions $\bar{A}_j$ which satisfy

\begin{equation}\label{redu}\partial_\omega \Psi_j =\bar{A}_j\Psi_j -\Psi_j A_j.
\end{equation}
%







\bigskip
$\bullet$ \textit{Reducibility:} If there exists $J\geq 1$ such that $\Psi_{j+1}=\Psi_j$ for all $j\geq J$, then, taking the limit in \eqref{redu} in $C^{k'}$, one finds a matrix $A_\infty$ satisfying

$$\partial_\omega \Psi_J= \bar{A}_\infty \Psi_J-\Psi_J A_\infty$$

so $\bar{A}_\infty$ is reducible, and therefore $A+F$ is reducible. $\Box$
%

%

\section{Application to the fibered rotation number}

In this section, we focus on the case $G=SL(2,\mathbb{R})$, which includes the very important example of Schrödinger cocycles.

\begin{prop} \label{redppck}Let $\bar{A}\in C^0(\mathbb{T}^d,sl(2,\mathbb{R}))$. There exists $k,D'$ only depending on $d,\kappa,\tau,\hat{\bar{A}}(0)$ and $\epsilon_0$ only depending on $d,\kappa,\tau,\hat{\bar{A}}(0),k$ such that if
\begin{itemize}
\item $\bar{A}\in C^k(\mathbb{T}^d,sl(2,\mathbb{R}))$, 
\item $\rho(\bar{A})\in DC_\omega(\tau)\cup \mathcal{M}_\omega$ 
\item and $\mid \mid \bar{A}-\hat{\bar{A}}(0)\mid \mid _k\leq \epsilon_0$, 
\end{itemize}
then the cocycle associated to $\bar{A}$ is reducible in $C^{k-D'}(\mathbb{T}^d,SL(2,\mathbb{R}))$.
\end{prop}

\dem We shall apply Theorem \ref{PRck} with $A=\hat{\bar{A}}(0)$ and $F=\bar{A}-\hat{\bar{A}}(0)$. Let $k_0,D'$ only depending on $d,\kappa,\tau,\hat{\bar{A}}(0)$  as in Theorem \ref{PRck}, $k\geq k_0$ and 
$\epsilon_0$ only depending on $d,\kappa,\tau,\hat{\bar{A}}(0),k$ as in Theorem \ref{PRck}. If $\mid \mid \bar{A}-\hat{\bar{A}}(0)\mid \mid _k\leq \epsilon_0$, there exists 
\begin{itemize}
\item $Z_\infty\in C^{k-D'}(\mathbb{T}^d,SL(2,\mathbb{R}))$,
\item $\bar{A}_\infty\in C^{k-D'}(\mathbb{T}^d,sl(2,\mathbb{R}))$
\end{itemize} 
such that 

$$\partial_\omega Z_\infty= \bar{A} Z_\infty-Z_\infty \bar{A}_\infty$$

\noindent
and $\bar{A}_\infty$ is the limit of a sequence of maps $(\bar{A}_j)$ which are reducible to $A_j$ by $\Psi_j$. 
%





\bigskip
\noindent
Let $(M_j), (\kappa_j), (N_j)$ be sequences as in Theorem \ref{PRck}. 
We shall proceed by contradiction; suppose that the cocycle associated to $\bar{A}$ is not reducible: then there exists a sequence $(j_l)_{l\geq 1}$  
such that for all $l$, $\Psi_{j_l+1}\neq \Psi_{j_l}$ (i.e $M_{j_l}\neq 0$).
Now by definition of the sequence $(M_j)$, for all $j$,

$$\rho(A_j)+2\pi \sum_{j'=1}^{j-1} \langle M_{j'},\omega\rangle=\rho(\bar{A}_j)$$ 

\noindent so we have

$$\mid \rho(\bar{A}_j)-\rho(\bar{A}_{j+1})\mid =\mid \rho(A_{j+1})-(\rho(A_j)-2\pi \langle M_j,\omega\rangle)\leq \kappa_j$$

\noindent and therefore, by continuity of the rotation number,

$$\mid \rho(\bar{A}_j)-\rho(\bar{A}_\infty)\mid \leq \sum_{j'\geq j}\kappa_{j'}.$$

\noindent 
Suppose $\rho(\bar{A})\in DC_\omega(\tau)$, then for all $m\in \mathbb{Z}^d\setminus\{0\}$, 

$$\mid \rho(\bar{A})-2\pi \langle m,\omega\rangle \mid \geq \frac{\kappa'}{\mid m\mid ^{\tau}}$$

\noindent
for some $\kappa'$; now there exists $M\in \mathbb{Z}^d$ such that for all $m\in \mathbb{Z}^d\setminus\{0\}$, 

$$\mid \rho(\bar{A}_\infty)-2\pi \langle m,\omega\rangle\mid 
=\mid \rho(\bar{A})-2\pi \langle m-M,\omega\rangle \mid.$$

\noindent
So,
for all $m\in \mathbb{Z}^d\setminus\{M\}$ and all $l$ big enough, 

\begin{equation}\label{infty}\mid \rho(A_{j_l})+2\pi \sum_{j'=1}^{j_l-1} \langle M_{j'},\omega\rangle -2\pi \langle m,\omega\rangle\mid +\sum_{j'\geq j_l}\kappa_{j'}
\geq \frac{\kappa'}{\mid m-M\mid^{\tau}}.
\end{equation}

\noindent
Now by definition of the sequence $(j_l)$, for all $l$, 

$$\mid \rho(A_{j_l})-2\pi \langle M_{j_l},\omega\rangle \mid 
<\kappa_{j_l}$$ 

\noindent
therefore

\begin{equation}\mid \rho(A_{j_l})+2\pi \sum_{j'=1}^{j_l-1} \langle M_{j'},\omega\rangle -2\pi \langle \sum_{j'=1}^{j_l} M_{j'},\omega\rangle\mid+\sum_{j'\geq j_l}\kappa_{j'}
 < 2\sum_{j'\geq j_l}\kappa_{j'}.
\end{equation}

\noindent
If one lets

\begin{equation}
\begin{split}\kappa'_l= 2\sum_{j'\geq j_l}\kappa_{j'}[\mid \sum_{j'=1}^{j_l} M_{j'}\mid+\mid M\mid ] ^{\tau}
\end{split}\end{equation}

\noindent then

\begin{equation}\kappa'_l> 2\sum_{j'\geq j_l}\kappa_{j'}(j_lR_{j_l-1}N_{j_l-1}+\mid M\mid )^{\tau}>0
\end{equation}

\noindent and 

\begin{equation}\label{decr}\mid \rho(A_{j_l})+2\pi \sum_{j'=1}^{j_l-1} \langle M_{j'},\omega\rangle+\sum_{j'\geq j_l}\kappa_{j_l} -2\pi \langle \sum_{j'=1}^{j_l} M_{j'},\omega\rangle\mid<\frac{\kappa'_l}{\mid \sum_{j'=1}^{j_l} M_{j'}-M\mid ^{\tau}}.
\end{equation}

\noindent The sequence $\kappa'_l$ also satisfies

$$ \kappa'_l\leq 2\sum_{j'\geq j_l}\kappa_{j'}[\sum_{j'=1}^{j_l} N_{j'}+\mid M\mid ]^{\tau}\leq c(j_lN_{j_l}+\mid M\mid )^{\tau}\kappa_{j_l}\leq c'(j_l\frac{1}{R_{j_l}^{2}}+\mid M\mid \kappa_{j_l}^{\frac{1}{\tau}})^\tau$$

\noindent where $c,c'$ do not depend on $l$. 
For $l$ big enough, \eqref{decr} contradicts \eqref{infty} since $(M_j)_{j\geq 1}$ is unbounded and $\kappa'_l$ tends to 0. 

\bigskip
\noindent In the case when $\rho(\bar{A})$ is rational with respect to $\omega$, since $\bar{A}$ and $\bar{A}_\infty$ are conjugated, $\rho(\bar{A}_\infty)$ is also rational with respect to $\omega$. Therefore, \eqref{infty} still holds, with $\kappa'=\kappa$, for $M$ such that $\rho(\bar{A}_\infty)=\langle M,\omega\rangle$, so one is led to the same contradiction. 
$\Box$

\bigskip
\noindent
Let us now consider the cocycle associated to the Schrödinger equation

\begin{equation}\label{schr}\frac{d}{dt}X(t)=(A_\lambda+F(t\omega)) X(t)
\end{equation}

\noindent 
where $A_\lambda=\left(\begin{array}{cc}
0 & -\lambda\\
1& 0\\
\end{array}\right)$ and $F(t\omega)=\left(\begin{array}{cc}
0 & V(t\omega)\\
0& 0\\
\end{array}\right)$ with $V\in C^k(\mathbb{T}^d)$ with $k$ to be determined later on.

\begin{thm}There exists $k_0$ only depending on $d,\kappa,\tau$ such that if $k\geq k_0$ and if $V\in C^k(\mathbb{T}^d)$, there exists $\epsilon_0$ only depending on $d,\kappa,\tau, k$ such that if $\mid \mid V\mid \mid _k \leq \epsilon_0$, then the cocycle which is solution of \eqref{schr} is
\begin{itemize}
\item almost reducible for all $\lambda$,
\item reducible for all $\lambda$ such that $\rho(A_\lambda+F)\in DC_\omega(\tau)\cup \mathcal{M}_\omega$. 
\end{itemize}
\end{thm}

\dem $\bullet$ First case: $\lambda\in [-2,2]$. The norm of $A_\lambda$ is then bounded independently of $\lambda$ so it is enough to apply Theorem \ref{PRck} with $A=A_\lambda$ and $F$ as above to deduce almost reducibility; to infer reducibility if $\rho(A_\lambda+F)\in DC_\omega(\tau)\cup \mathcal{M}_\omega$, apply Proposition \ref{redppck} with $\bar{A}=A_\lambda+F$.

\bigskip
\noindent
$\bullet$ Second case: $\mid \lambda\mid >2$. Letting $Y(t)=\left(\begin{array}{cc}
\frac{1}{2} & -\frac{\sqrt{\lambda}}{2}\\
\frac{1}{2} & \frac{\sqrt{\lambda}}{2}\\
\end{array}\right)X(t)$, one has

\begin{equation}
\begin{split} Y'(t)&=\left(\begin{array}{cc}
\frac{1}{2} & -\frac{\sqrt{\lambda}}{2}\\
\frac{1}{2} & \frac{\sqrt{\lambda}}{2}\\
\end{array}\right)
\left(\begin{array}{cc}
0 & V(t\omega)-\lambda\\
1& 0\\
\end{array}\right)\left(\begin{array}{cc}
1 & 1\\
-\frac{1}{\sqrt{\lambda}} & \frac{1}{\sqrt{\lambda}}\\
\end{array}\right)Y(t)\\
&=(\tilde{A}(\lambda)+\tilde{F}(\lambda,t\omega))Y(t)
\end{split}
\end{equation}

\noindent 
with 

$$\tilde{A}(\lambda)
=\left(\begin{array}{cc}
0 & -\sqrt{\lambda} \\
\sqrt{\lambda} & 0\\
\end{array}\right)$$ 

\noindent
and 

$$\tilde{F}(\lambda,t\omega)=\left(\begin{array}{cc}
-\frac{V(t\omega)}{2\sqrt{\lambda}} & \frac{V(t\omega)}{2\sqrt{\lambda}}\\
-\frac{V(t\omega)}{2\sqrt{\lambda}} & \frac{V(t\omega)}{2\sqrt{\lambda}}\\
\end{array}\right).$$

\noindent
Thus, one can apply Theorem \ref{PRck} with $A=\tilde{A}(\lambda)$ and $F(t\omega)=\tilde{F}(\lambda,t\omega)$ to get almost reducibility if $V$ is bounded in the $C^k$ topology by some constant depending only on $d,\kappa,\tau,k$. 
One can also apply Proposition \ref{redppck} with $\bar{A}=\tilde{A}(\lambda)+\tilde{F}(\lambda, t\omega)$ to get reducibility in the case when $\rho(A_\lambda+F)\in DC_\omega(\tau)\cup \mathcal{M}_\omega$, since $\rho(A_\lambda+F)=\rho (\tilde{A}(\lambda)+\tilde{F}(\lambda,.))$. 
$\Box$

\section{Acknowledgments} The author would like to thanks the Institut Mittag-Leffler in Stockholm for its hospitality, as well as L.H. Eliasson for useful discussions.

\end{document}